\documentclass[11pt,abstract]{scrartcl}
\usepackage{amsmath,amssymb,amsthm}
\usepackage{physics}
\usepackage[final]{microtype}
\usepackage{multirow}
\usepackage{xcolor}
\usepackage{tikz}
\usetikzlibrary{arrows.meta, external, calc, spy, pgfplots.groupplots}
\tikzset{external/only named=true}
\usepackage{pgfplots}
\pgfplotsset{compat=newest}

\usepackage[backend=biber,hyperref,backref,bibencoding=utf8,defernumbers=true,style=numeric-comp,giveninits=true,doi=true,url=false,isbn=false]{biblatex}
\renewbibmacro{in:}{\ifentrytype{article}{}{\printtext{\bibstring{in}\intitlepunct}}}
\usepackage{csquotes}
\usepackage{hyperref}
\usepackage{cleveref}

\renewcommand{\vec}[1]{{\boldsymbol{#1}}}

\newcommand{\abssec}[1]{\noindent\small {\bfseries #1\quad}\ignorespaces}
\renewenvironment{abstract}{\abssec{Abstract}}{\par\vspace{1em}}

\makeatletter
\let\blx@rerun@biber\relax
\makeatother

\definecolor{unianthrazit}{rgb}{0.416,0.408,0.435}
\definecolor{uniorange}{rgb}{0.929,0.431,0.0}
\definecolor{adjyellow}{rgb}{0.929,0.651,0.0}
\definecolor{adjred}{rgb}{0.898,0.0,0.086}
\definecolor{compblue}{rgb}{0.004,0.537,0.561}
\pgfplotsset{colormap={unibwthesis}{color(0cm)=(compblue);color(1cm)=(adjyellow);color(2cm)=(uniorange);color(3cm)=(adjred)}, colormap name=unibwthesis}

\title{Pressure-robustness for the Stokes equations on anisotropic meshes}
\date{\today}
\author{Volker Kempf}

\begin{document}
\maketitle

\begin{abstract}
	Pressure-robustness has been widely studied since the conception of the notion  and the introduction of the reconstruction approach for classical mixed methods in \cite{Linke2014}.
	Using discretizations capable of yielding velocity solutions that are independent of the pressure approximation has been recognized as essential, and a large number of recent articles attest to this fact, e.g., \cite{LinkeMerdonNeilan2020,AinsworthParker2021}.
	Apart from the pressure-robustness aspect, incompressible flows exhibit anisotropic phenomena in the solutions which can be dealt with by using anisotropic mesh grading.
	The recent publications \cite{ApelKempf2021,ApelKempfLinkeMerdon2021} deal with the combination of both challenges.
	We briefly revisit the results from \cite{ApelKempfLinkeMerdon2021} and provide an insightful new numerical example.
\end{abstract}

\section{Stokes equations, modified Crouzeix--Raviart discretization, error estimates}
We consider the Stokes equations in a domain $\Omega$ with a kinematic viscosity $\nu> 0$ and homogeneous Dirichlet boundary conditions. The weak form of this problem is: Find $(\vec{u},p)\in\vec{H}_0^1(\Omega)\times L_0^2(\Omega)$ so that
\begin{align*}
&\nu (\nabla \vec{u}, \nabla\vec{v}) - (\nabla\cdot\vec{v},p) - (\nabla\cdot\vec{u},q) = (\vec{f},\vec{v}) &&\forall (\vec{v}, q) \in \vec{H}_0^1(\Omega)\times L_0^2(\Omega).
\end{align*}
By using the Helmholtz decomposition of the data $\vec{f}=\mathbb{P} \vec{f}+\nabla\phi$ and looking at the problem in the subspace of divergence free functions, it can be seen that the velocity solution is independent of the gradient part $\nabla \phi$ of the data, see \cite{Linke2014}.
Most classical mixed methods do not mimic this in the discrete setting.
Also, solutions of incompressible flow problems exhibit anisotropic structures like boundary layers or edge singularities which can be dealt with by using anisotropically graded meshes, see \cite{ApelKempfLinkeMerdon2021,ApelKempf2021}.

The modified Crouzeix--Raviart method was introduced in \cite{Linke2014}, and extended to anisotropic meshes in \cite{ApelKempfLinkeMerdon2021,ApelKempf2021}. The idea behind the modification is to alter the discrete test functions in the linear form in order to recover the $L^2$-orthogonality of the test function and the gradient part of the Helmholtz-decomposition of the data function. This leads to the discretization
\begin{align*}
&\nu (\nabla_h \vec{u}_h, \nabla_h\vec{v}_h) - (\nabla_h\cdot\vec{v}_h,p_h) - (\nabla_h\cdot\vec{u}_h,q_h) = (\vec{f},I_h \vec{v}_h) &&\forall (\vec{v}_h, q_h) \in \vec{X}_h\times Q_h,
\end{align*}
where $\vec{X}_h$ is the space of piece-wise linear functions that are continuous only at the barycenter of element interfaces, $Q_h$ contains the piece-wise constant functions and $\nabla_h$, $\nabla_h\cdot$ are the broken gradient and divergence. 
With $I_h=\operatorname{id}$ we get the standard Crouzeix--Raviart method (CR), while for the pressure-robust modification we can choose $I_h$ as the lowest-order Raviart--Thomas (CR-RT) or Brezzi--Douglas--Marini (CR-BDM) interpolation operators, see \cite{Linke2014} and for anisotropic meshes \cite{ApelKempf2021,ApelKempfLinkeMerdon2021,ApelKempf2020}. 
Assuming for the solution $(\vec{u},p)\in \vec{H}^2(\Omega)\times H^1(\Omega)$ and a mesh that satisfies the maximum angle condition, i.e., the maximum angle in the mesh is uniformly bounded away from $\pi$, the error analysis for the modified methods yields the pressure-independent velocity estimate, see \cite[Theorem 4.1]{ApelKempfLinkeMerdon2021},
\begin{align*}
\norm{\vec{u}-\vec{u}_h}_{1,h} &\leq 2 \inf_{\vec{v}_h\in\vec{X}^0_h} \norm{\vec{u} - \vec{v}_h}_{1,h} + Ch\abs{\vec{u}}_{\vec{H}^2(\Omega)},\text{ where } \norm{\vec{v}_h}^2_{1,h} = (\nabla_h\vec{v}_h, \nabla_h\vec{v}_h)\text{ and }
\end{align*}
$\vec{X}_h^0$ is the space of discretely divergence free functions from $\vec{X}_h$.
The estimate for the pressure is similar to the CR method.

\section{Numerical example}
\begin{figure}[tbp]
	\centering
	\includegraphics{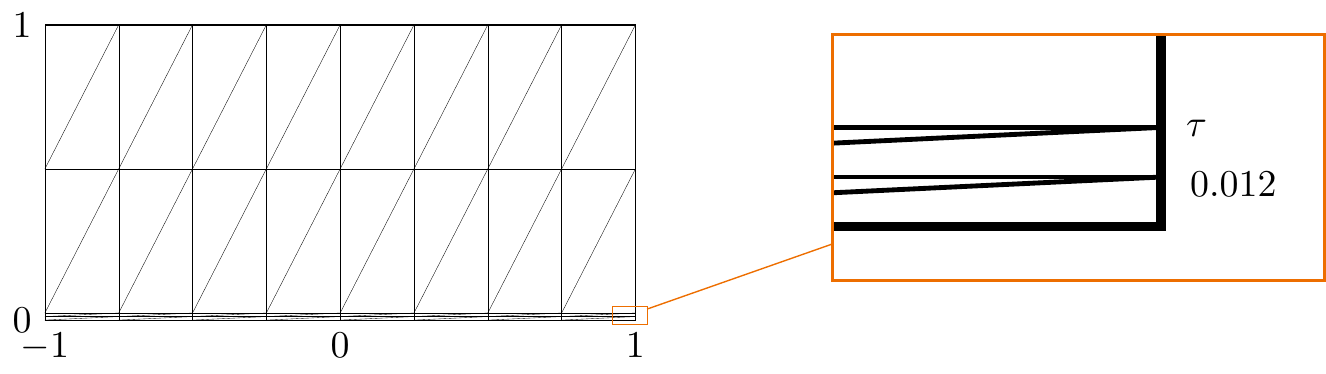}
	\caption{Shishkin-type mesh with transition point $\tau$ set to layer width $\delta$, resulting for $\nu=10^{-4}$ in an aspect ratio of approx. $21.4$.}
	\label{fig:shishkin_mesh}		
\end{figure}

\begin{figure}[htbp]
	\centering
	\includegraphics{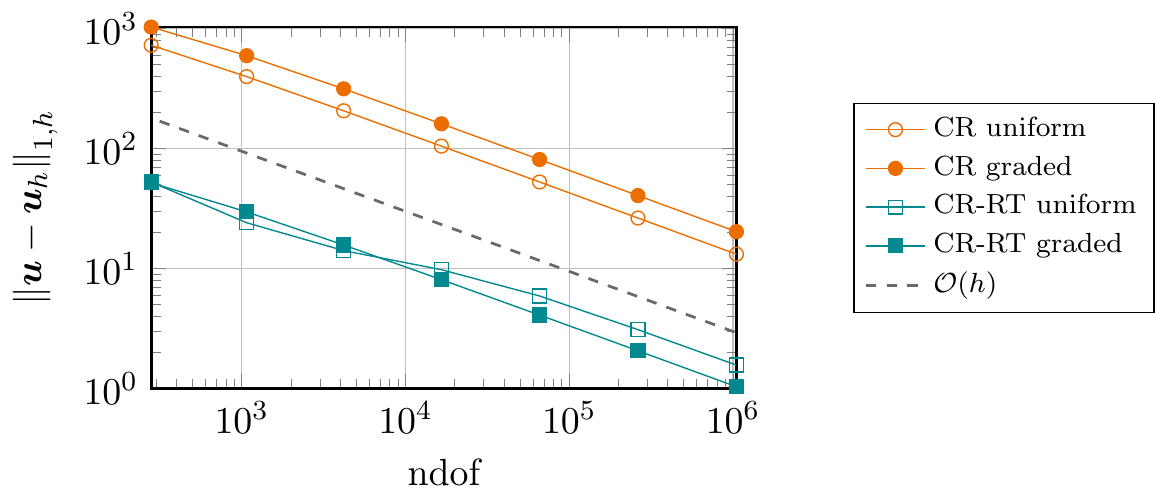}\hfill
	\includegraphics{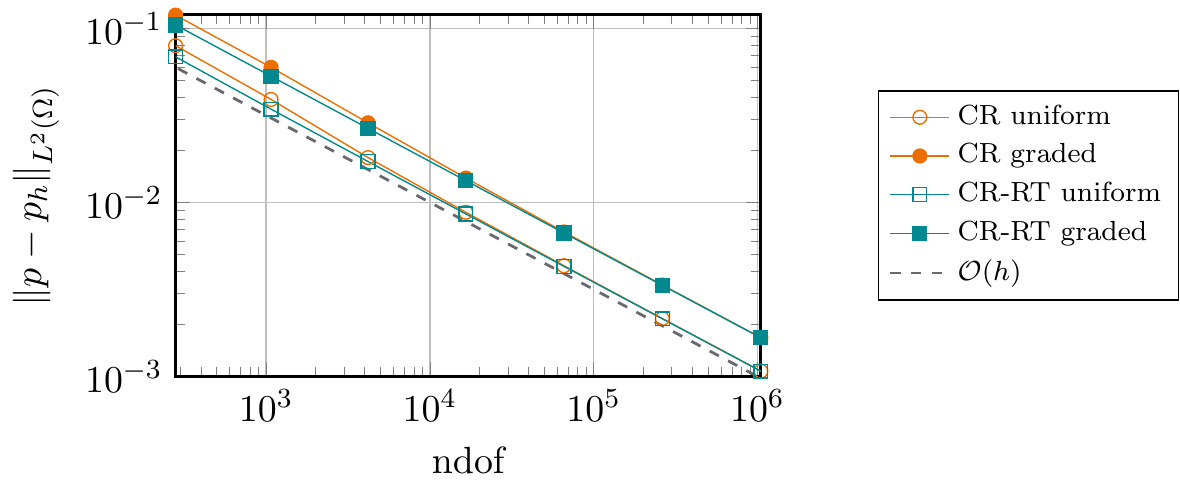}
	\caption{Convergence plots for the standard and modified Crouzeix--Raviart methods.}
	\label{fig:stag_point_convergence}		
\end{figure}

\begin{figure}[tbp]
	\centering
			\includegraphics{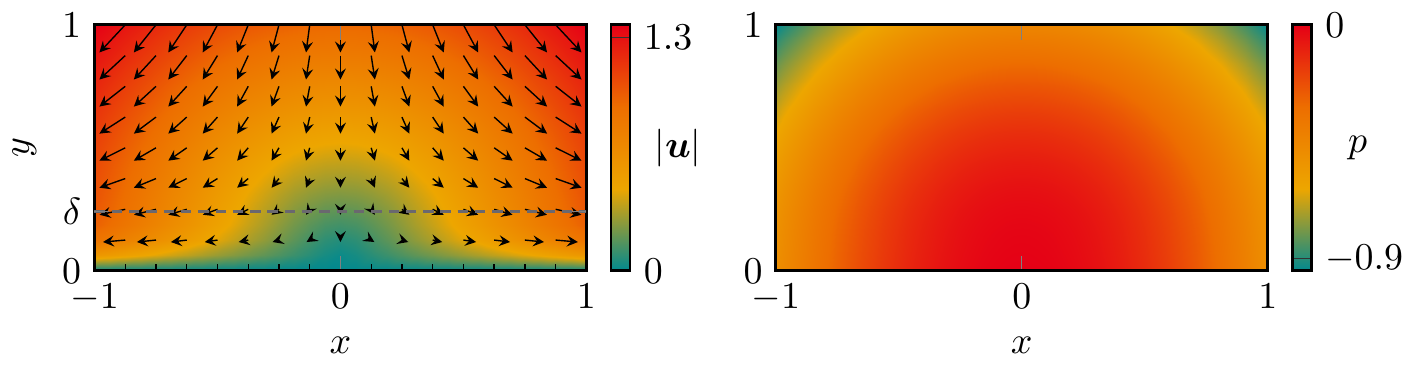}%
	\caption{Velocity and pressure solution for $\nu=0.01$, $a=1$, $P_0=0$.}
	\label{fig:exact_sol}
\end{figure}

\begin{figure}[tbp]	
		\centering
		\includegraphics{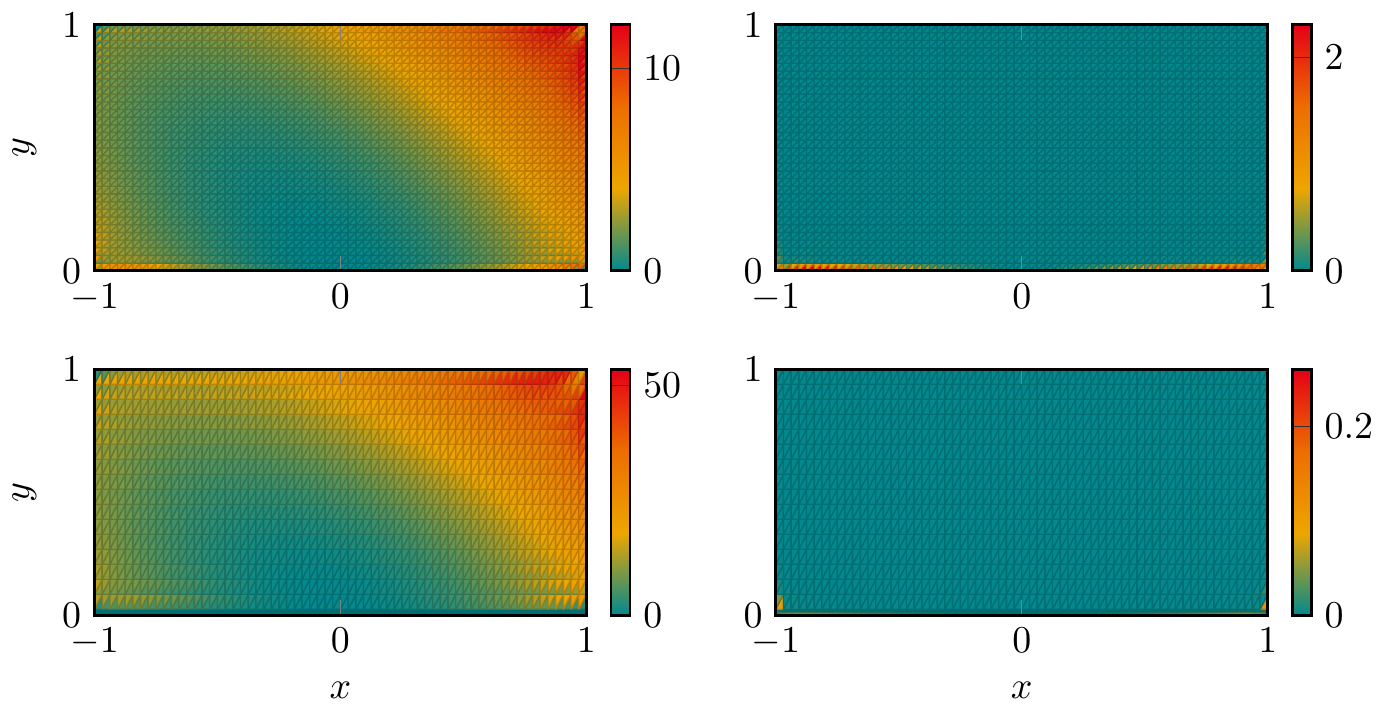}%
		\caption{Element-wise error $\norm{\nabla(\vec{u}-\vec{u}_h)}_{\vec{L}^2(T)}$ for CR (left), CR-RT (right) methods and uniform (top), graded (bottom) meshes.}
		\label{fig:elementwise_error}	
\end{figure}

We investigate a stagnation point flow in two dimensions, see \cite{SchlichtingGersten2006}, i.e., a stationary laminar flow against a wall.
Using the domain $\Omega=(-1,1)\times(0,1)$ and parameters $a\in \mathbb{R}_+$, $P_0\in\mathbb{R}$, an exact solution of the Navier--Stokes equations is given by
\begin{equation*}
\vec{u}(x,y) = \begin{pmatrix}axf'(\eta) \\ -\sqrt{a\nu}f(\eta)\end{pmatrix},\  p(x,y)=P_0-\frac{a^2}{2}\left(x^2 + \frac{2\nu}{a} \left(f'(\eta)+\frac{1}{2}f^2(\eta)\right)\right),\ \eta = \sqrt{\frac{a}{\nu}}y.
\end{equation*}
The function $f$ is defined via the boundary value problem $f''' + ff'' + 1 - f'^2=0$, $f(0)=f'(0)=0$, $\lim_{\eta\to\infty} f'(\eta) = 1$.
This BVP can be solved to high precision, e.g. using a shooting method, which yields a numerically highly accurate solution, see Fig.~\ref{fig:exact_sol}, that can be used for the convergence study.	
The velocity solution shows a boundary layer structure near $y=0$, with a width of $\delta=2.4\sqrt{\nu a^{-1}}$, which means that anisotropic meshes should be useful.
In order to stay in the Stokes framework, we use the exact solution to define a data function $\vec{f}=-(\vec{u}\cdot\nabla)\vec{u}$ to examine the performance of the discretizations.
For the parameter choices $\nu=10^{-4}$, $a=1$, $P_0=0$ we compare the numerical results of the CR and CR-RT methods on uniform structured meshes and meshes of Shishkin-type, as pictured in Fig.~\ref{fig:shishkin_mesh}. 
In Fig.~\ref{fig:stag_point_convergence} the corresponding convergence plots are shown.
The modified method does not show an advantage for the pressure approximation and the graded meshes increase the pressure error, which is to be expected, given that the pressure does not have a boundary layer, see Fig.~\ref{fig:exact_sol}.
For the velocity on the other hand, the pressure-robustness of the CR-RT method reduces the error by about one order of magnitude.
On graded meshes there is additional improvement of the velocity solution with the CR-RT method, while the error increases for the CR method, which at first seems unreasonable.
To find a reason we have a look at the element-wise error in Fig.~\ref{fig:elementwise_error}.
The left hand side of the figure shows that the error of the CR method is distributed throughout the domain, so that the grading can not work as intended, as seen in the bottom left plot.
On the other side of the figure we see that the error of the CR-RT method is concentrated in the boundary layer, so the grading does indeed improve the solution.

The numerical results of this example show that using a pressure-robust method in combination with anisotropically graded meshes holds great potential and can be highly advantageous over standard discretization schemes.

\printbibliography

\end{document}